%% file: HET.tex
\documentclass[12pt]{article}
\usepackage{a4wide,graphics,bm,color}
\usepackage{amsmath,latexsym,amssymb,amsfonts}
\usepackage{ijbc}
\usepackage{tabls}
\usepackage{array,delarray}
\usepackage{graphicx}
\begin{document}
\begin{titlepage}
  \title{CONTINUATION OF CONNECTING ORBITS IN 3D-ODES: (I)
    POINT-TO-CYCLE CONNECTIONS} \author{E.J. DOEDEL$^1$, B.W.
    KOOI$^2$, YU.A. KUZNETSOV$^3$, and G.A.K. van VOORN$^2$
    \\ \\
    {\it $^1$Department of Computer Science, Concordia University,}\\
    {\it 1455 Boulevard de Maisonneuve O., Montreal, Quebec, H3G 1M8, Canada}\\
    {\tt doedel@cs.concordia.ca}
    \\[1mm]
    {\it $^2$Department of Theoretical Biology,  Vrije Universiteit,}\\
    {\it de Boelelaan 1087, 1081 HV Am\-ster\-dam,
      the Netherlands}\\
    {\tt kooi@bio.vu.nl, george.van.voorn@falw.vu.nl}
    \\[1mm]
    {\it $^3$Department of Mathematics, U\-trecht University}\\
    {\it Budapestlaan 6, 3584 CD Utrecht, the Netherlands}\\
    {\tt kuznet@math.uu.nl}
    \\[1mm]
    } \date{\today} \maketitle \abstract{We propose new methods for
    the numerical continuation of point-to-cycle connecting orbits in
    3-dimensional autonomous ODE's using projection boundary
    conditions.  In our approach, the projection boundary conditions
    near the cycle are formulated using an eigenfunction of the
    associated adjoint variational equation, avoiding costly and
    numerically unstable computations of the monodromy matrix.
    The equations for the eigenfunction
    are included in the defining boundary-value problem, allowing
    a straightforward implementation in {\sc auto}, in which only 
    the standard features of the software are
    employed. Homotopy methods to find connecting orbits are discussed
    in general and illustrated with several examples, including the
    Lorenz equations. Complete {\sc auto} demos, which can be
    easily adapted to any autonomous 3-dimensional ODE system, 
    are freely a\-vai\-la\-ble.}
\end{titlepage}

\twocolumn
\include{hetero}

\bibliographystyle{ijbc} 
\clearpage
\bibliography{refs}

\include{appendix}

\end{document}

%% file: hetero.tex
\section{Introduction}\label{sec:intro}

Many interesting phenomena in ODE systems can only be understood by
analyzing global bifurcations. Examples of such are the occurrence and
disappearance of chaotic behaviour. For example, the classical Lorenz
attractor appears in a sequence of bifurcations, where {\em homoclinic
orbits} connecting a saddle equilibrium to itself and {\em
heteroclinic orbits} connecting an equilibrium point with a saddle
cycle, are involved (Afraimovich et al., 1977\nocite{AfrBykShil1977}). 
In the ecological context, Boer et al.
(1999\nocite{Boeretal1999}, 2001\nocite{Boeretal2001}) showed that
regions of chaotic behaviour in parameter space in some food chain
models are bounded by bifurcations of point-to-cycle and
cycle-to-cycle connections.

Thus, in order to gain more knowledge about the global bifurcation
structure of a model, information is required on the existence of
\textit{homoclinic} and \textit{heteroclinic} connections between
equilibria and/or periodic cycles. The first type is a connection that
links an equilibrium or a cycle to itself (asymptotically bi-stable, so
it necessarily has nontrivial stable and unstable invariant manifolds). 
The second type is a connection
that links an equilibrium or a cycle to another equilibrium or cycle.

The continuation of connecting orbits in ODE systems has been
notoriously difficult. Doe\-del and Friedman
(1989\nocite{DoedelFriedman1989}) and Beyn (1990\nocite{Beyn1990})
developed direct numerical methods for the computation of orbits
connecting equilibrium points and their associated parameter values,
based on truncated boundary value problems with {\em projection
boundary conditions}.  Moreover, Doe\-del, Fried\-man and Mon\-tei\-ro
(1993\nocite{DoFrMo1994}) have proposed efficient methods to find
starting solutions by successive continuations (homotopies).  These
continuation methods have been implemented in {\sc HomCont}, as
incorporated in {\sc auto} (Doedel et al.,
1997\nocite{Doedeletal1997}; Cha\-mpneys and Kuznetsov,
1994\nocite{ChampneysKuznetsov1994}; Cha\-mpneys et al.,
1996\nocite{Champneysetal1996}). {\sc HomCont} is only suitable for
the continuation of homoclinic point-to-point and heteroclinic
point-to-point connections.

More recently, significant progress has been made in the continuation
of homoclinic and heteroclinic connections involving cycles.  Die\-ci
and Rebaza (2004\nocite{DieciRebaza2004a,DieciRebaza2004}) developed a
method based on earlier works by Beyn (1994\nocite{Beyn1994}) and
Pampel (2001\nocite{Pampel2001}).  Their method is also based on
projection boundary conditions, but uses an {\em ad hoc} multiple
shooting technique and requires the numerical determination of the
monodromy matrix associated with the periodic cycles involved in the
connection.

In this paper, we propose new methods for the numerical continuation
of point-to-cycle connections in 3-dimensional autonomous ODE's using
projection boundary conditions.  In our approach, the projection
boundary conditions near each cycle are formulated using an
eigenfunction of the associated adjoint variational equation, avoiding
costly and numerically unstable computation of the monodromy matrix.
Instead, the equations for the eigenfunction
are included in the defining boundary-value problem, allowing a
straightforward implementation in {\sc auto}. 

This paper is organized as follows. In Section \ref{sec:algorithm} we
recall basic properties of the projection boundary condition method
to continue point-to-cycle connections. In Section
\ref{heteroclinicBVP} this method is adapted to efficient
numerical implementation in a special -- but important -- 3D case.
Homotopy methods to find connecting orbits are discussed in Section
\ref{sec:homotopy}. Section \ref{sec:implementation} demonstrates that
the algorithms allow for a straightforward implementation in
{\sc auto}, using only the basic features of this software.  
Three well-known examples (the three-dimensional Lorenz
system, the electronic circuit model of Freire et al.,
1993\nocite{Freireetal1993}, and the standard three-level food chain 
model based on the Rosen\-zweig-MacArthur (1963\nocite{RosenzweigMacArthur1963})
system) are used in Section \ref{sec:results} to illustrate the power 
of the new methods.

This is Part I of a sequel of two papers. Part II will deal with
cycle-to-cycle connections in 3D systems.

\section{Truncated BVP's with projection BC's}\label{sec:algorithm}
Before presenting a BVP for a point-to-cycle  
connection, we set up some notation.  
\begin{figure*}[htbp]
\begin{center}
  \includegraphics[width=17.0cm]{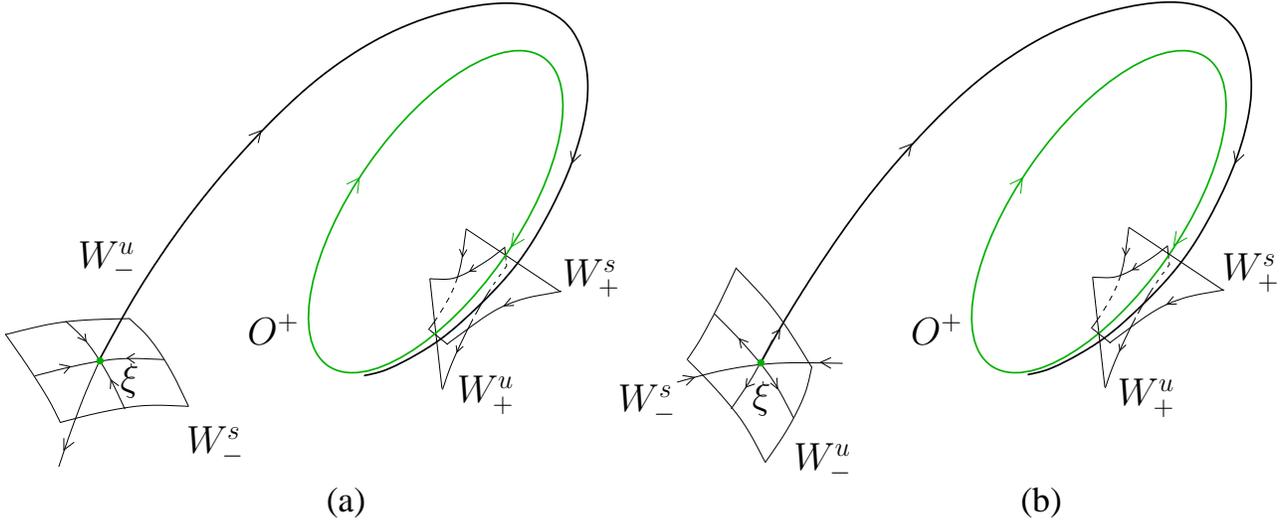}
\end{center}
\caption{Point-to-cycle connecting orbits in ${\mathbb R}^3$: (a)
$n^{-}_u=1$; (b) $n^{-}_u=2$.}
\label{fig:connections} 
\end{figure*}

Consider a general system of ODE's
\begin{equation}
\label{eqn:ODE}
  \frac{du}{dt} = f(u,\alpha),
\end{equation}
where $f: \mathbb{R}^n \times \mathbb{R}^p \to \mathbb{R}^n$ is a
sufficiently smooth function of the state variables $u\in
\mathbb{R}^n$ and the control parameters $\alpha \in
\mathbb{R}^p$. Denote by $\varphi^t$ the (local) flow generated by
(\ref{eqn:ODE})\footnote{Whenever possible, we will not indicate
explicitly the dependence of various objects on system parameters.}.

Let $O^{-}$ be either a saddle or a saddle-focus equilibrium, say
$\xi$, and let $O^{+}$ be a hyperbolic saddle limit cycle of
(\ref{eqn:ODE}).  A solution $u(t)$ of (\ref{eqn:ODE}) defines a {\em
connecting orbit} from $O^{-}$ to $O^{+}$ if
\begin{equation}
\label{eqn:ASSYMPT}
\lim_{t \to \pm \infty}{\rm dist}(u(t), O^{\pm}) = 0
\end{equation}
(see Figure \ref{fig:connections} for illustrations).  Since
$u(t+\tau)$ satisfies (\ref{eqn:ODE}) and (\ref{eqn:ASSYMPT}) for any
phase shift $\tau$, an additional scalar phase condition
\begin{equation}
\label{eqn:PHASE}
\psi[u,\alpha] = 0
\end{equation}
is needed to ensure uniqueness of the connecting orbit.
This condition will be specified later.

For numerical approximation, the asymptotic conditions
(\ref{eqn:ASSYMPT}) are replaced by {\em projection boundary
conditions} at the end-points of a large {\em truncation interval}
$[\tau_{-},\tau_{+}]$: The points $u(\tau_{-})$ and $u(\tau_{+})$ are
required to belong to the linear subspaces that are tangent to the
unstable and stable invariant manifolds of $O^{-}$ and $O^{+}$,
respectively.

Let $n^{-}_u$ be the dimension of the unstable invariant manifold
$W^u_{-}$ of $\xi$, {\it i.e.}, the number of eigenvalues
$\lambda^{-}_u$ of the Jacobian matrix $f_u=D_uf$ evaluated at the
equilibrium which satisfy
$$
\Re(\lambda^{-}) >0.
$$

Denote by $x^{+}(t)$ a periodic solution (with minimal period $T^{+}$)
corresponding to $O^{+}$ and introduce the {\em monodromy matrix}
$$
M^{+}=\left.D_x \varphi^{T^{+}}(x)\right|_{x=x^{+}(0)}, 
$$
{\it i.e.}, the linearization matrix of the $T^{+}$-shift along orbits of
(\ref{eqn:ODE}) at point $x_0^{+}=x^{+}(0) \in O^{+}$. Its eigenvalues
$\mu^{+}$ are called the {\em Floquet multipliers}; exactly one of
them equals 1, due to the assumption of hyperbolicity.  
Let $m_s^{+}=n_s^{+}+1$ be the
dimension of the stable invariant manifold $W^s_{+}$ of the cycle
$O^{+}$; here $n_s^{+}$ is the number of its multipliers satisfying
$$
|\mu^{+}| < 1.
$$ A necessary condition to have an isolated {\em family} of
point-to-cycle connecting orbits of (\ref{eqn:ODE}) is that (see Beyn
(1994\nocite{Beyn1994}))
\begin{equation}
\label{eqn:FREEPARS1}
p=n-m_s^{+}-n_u^{-}+2
\end{equation}

The projection boundary conditions in this case can be written as
\begin{subequations}\label{eqn:BVPp2c}
\begin{align}
    L^{-} (u(\tau_{-}) - \xi) & = 0 \;, \label{eqn:BVPp2cunst}\\
    L^{+} (u(\tau_{+}) - x^{+}(0)) & = 0 \;, \label{eqn:BVPp2cst}
\end{align}
\end{subequations}
where $L^{-}$ is a $(n-n^{-}_u) \times n$ matrix whose rows form a
basis in the orthogonal complement of the linear subspace that is
tangent to $W^u_{-}$ at $\xi$. Similarly, $L^{+}$ is a $(n-m^{+}_s)
\times n$ matrix, such that its rows form a basis in the orthogonal
complement to the linear subspace that is tangent to $W^s_{+}$ of
$O^{+}$ at $x^{+}(0)$.

It can be proved that, generically, the truncated BVP composed of
(\ref{eqn:ODE}), a truncation of (\ref{eqn:PHASE}), and
(\ref{eqn:BVPp2c}) has a unique solution family
$(\hat{u},\hat{\alpha})$, provided that (\ref{eqn:ODE}) has a
connecting solution family satisfying (\ref{eqn:PHASE}) and
(\ref{eqn:FREEPARS1}).

The truncation to the finite interval $[\tau_{-},\tau_{+}]$ implies an
error. If $u$ is a generic connecting solution to (\ref{eqn:ODE}) at
parameter value $\alpha$, then the following estimate holds:
$$
\|(u|_{[\tau_{-},\tau_{+}]},\alpha) - (\hat{u},\hat{\alpha})\|\leq C {\rm e}^
{-2\min(\mu_{-}|\tau_{-}|,\mu_{+}|\tau_{+}|)},
$$ where $\|\cdot\|$ is an appropriate norm in the space
$C^1([\tau_{-},\tau_{+}],{\mathbb R}^n) \times {\mathbb R}^p$,
$u|_{[\tau_{-},\tau_{+}]}$ is the restriction of $u$ to the truncation
interval, and $\mu_{\pm}$ are determined by the eigenvalues of the
Jacobian matrix and the monodromy matrix. See Pampel
(2001\nocite{Pampel2001}) and Die\-ci and Rebaza
(2004\nocite{DieciRebaza2004a,DieciRebaza2004}) for exact
formulations, proofs, and references to earlier contributions.

\section{New defining systems in $\mathbb{R}^3$}
\label{heteroclinicBVP}
Here we explain how the projection boundary conditions
(\ref{eqn:BVPp2c}) can be implemented efficiently in a special -- but
important -- case $n=3$. Thereafter we specify the defining system
used to continue connecting orbits in 3D-ODE example systems with
\textsc{auto}. A saddle cycle $O^+$ in such systems always has
$m_s^{+}=m_u^{+}=2$.

\subsection{The equilibrium-related part}
The equilibrium point $\xi$, an appropriate solution of $f(\xi,\alpha)
= 0$, cannot be found by time-integration methods because it is a
saddle. There are two different types of saddle equilibria that can be
connected to saddle cycles in 3D-ODE's. These are distinguished by the
dimension $n^{-}_u$ of the unstable invariant manifold $W^u_{-}$ of
$\xi$: We have either $n^{-}_u=1$ or $n^{-}_u=2$ (see
Figure~\ref{fig:connections}). In the former case, the connection is
structurally unstable (has codim 1) and, according to
(\ref{eqn:FREEPARS1}), we need two free system parameters for its
continuation ($p=2$). In the latter case, however, the connection is
structurally stable and can be continued, generically, with one system
parameter ($p=1$). There is a small difference in the implementation
of the projection boundary condition (\ref{eqn:BVPp2cunst}) in these
two cases.
\begin{figure*}[ht]
\begin{center}
\includegraphics[width=16.0cm]{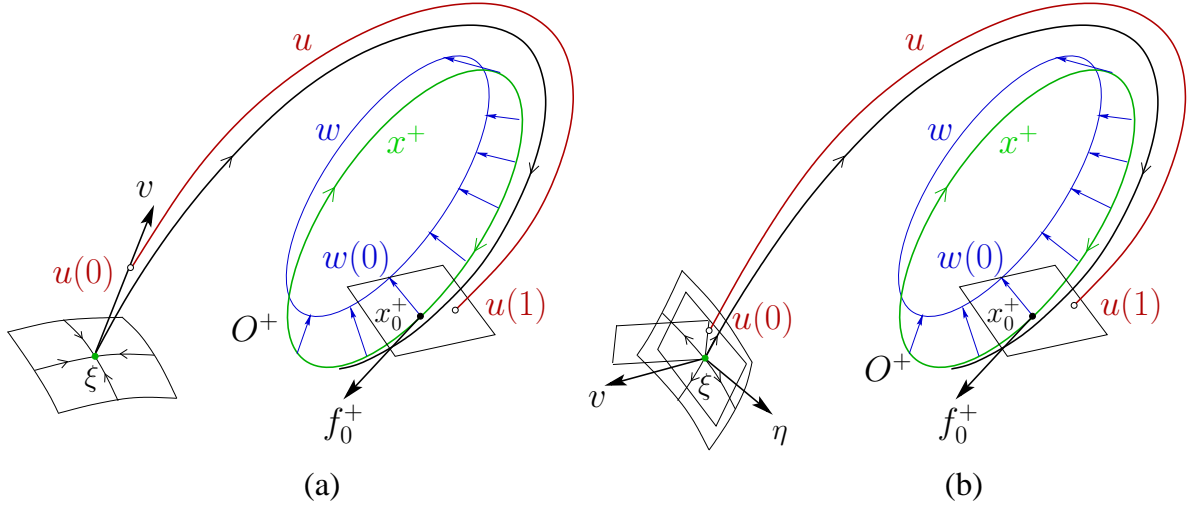}
\end{center}
\caption{BVP's to approximate connecting orbits: (a) $n^{-}_u=1$; (b)
$n^{-}_u=2$.}
\label{fig:conBVP} 
\end{figure*}

If $n^{-}_u=1$ (see Figure \ref{fig:conBVP}(a)), then the following
{\em explicit projection boundary condition} replaces
(\ref{eqn:BVPp2cunst}):
\begin{equation}
\label{eqn:PROJEXP1}
u(\tau_{-})=\xi + \varepsilon v,
\end{equation}
where $\varepsilon >0$ is a given small number, and $v \in {\mathbb
  R}^3$ is a unit vector that is tangent to $W^u_{-}$ at $\xi$. Notice
  that this fixes the phase of the connecting solution $u$, so that
  (\ref{eqn:PHASE}) becomes (\ref{eqn:BVPp2cunst}) in this case. The
  vector $v$ in (\ref{eqn:PROJEXP1}) is, of course, a normalized
  eigenvector associated with the unstable eigenvalue $\lambda_u >0$
  of the Jacobian matrix $f_u$ evaluated at the equilibrium. Hence, we
  can use the following algebraic system to continue $\xi,v$ and
  $\lambda_u$ simultaneously:
\begin{equation}\label{eqn:E1}
  \left\{
    \begin{array}{rcl}
      f(\xi,\alpha) & = &0\;,\\
      f_{\xi}(\xi,\alpha) v - \lambda_u v & = &0\;,\\
      \langle v,v \rangle - 1 & = &0 \;,
    \end{array}
  \right.
\end{equation}
where $\langle x,u \rangle = x^{\rm T} u$ is the standard scalar
product in $\mathbb{R}^n$.

If $n^{-}_u=2$ (see Figure \ref{fig:conBVP}(b)), then $W^{u}_{-}$ is
orthogonal to an eigenvector $v$ of the {\em transposed} Jacobian
matrix $f_u^{\rm T}$ corresponding to its eigenvalue $\lambda_s<0$, so
that (\ref{eqn:BVPp2cunst}) can be written as
\begin{equation}
\label{eqn:PROJEXP2}
\langle v, u(\tau_{-})-\xi \rangle = 0.
\end{equation}
To continue $\xi,v$, and $\lambda_s$, we use a system similar to
(\ref{eqn:E1}), namely:
\begin{equation}
\label{eqn:E2}
  \left\{
    \begin{array}{rcl}
      f(\xi,\alpha) & = &0\;,\\
      f^{\rm T}_{\xi}(\xi,\alpha) v - \lambda_s v & = &0\;,\\
      \langle v,v \rangle - 1 & = &0 \;.
    \end{array}
  \right.
\end{equation}
As a variant of the phase condition (\ref{eqn:PHASE}) in this case, we can use
the linear condition
\begin{equation}
\label{eqn:PROJEXP3}
\langle \eta, u(\tau_{-})-\xi \rangle = 0,
\end{equation}
which places the starting point of the truncated connecting solution
in a plane containing the equilibrium $\xi$ and orthogonal to a fixed
vector $\eta$ (not collinear with $v$). 

\subsection{The cycle and eigenfunctions}

The heteroclinic connection is linked on the other side to a saddle
limit cycle $O^+$ (see Figure \ref{fig:conBVP}). Thus, we also need a
BVP to compute it.  We use the standard periodic BVP:
\begin{equation}
\label{eqn:L}
    \left\{
      \begin{array}{ll}
        \dot{x}^+ - f(x^+,\alpha) & =  0 \;,\\
         x^+(0) - x^+(T^{+}) & =  0\;, 
      \end{array}
    \right.  
\end{equation}
which is augmented by an appropriate phase condition that makes its
solution unique. This phase condition is actually a boundary condition
for the truncated connecting solution, and will be introduced below.

To set up the projection boundary condition for the truncated
connecting solution $u$ near $O^+$, we also need a vector, say $w(0)$,
that is orthogonal at $x(0)$ to the stable manifold $W^s_{+}$ of the
saddle limit cycle $O^{+}$ (see Figure \ref{fig:conBVP}). It is well
known that $w(0)$ can be obtained from an {\em eigenfunction} $w(t)$
of the {\em adjoint variational problem} associated with
(\ref{eqn:L}), corresponding to its eigenvalue
$$
\mu=\frac{1}{\mu^{+}_u},
$$
where $\mu^{+}_u$ is a multiplier of the monodromy matrix $M^+$ satisfying
$$
|\mu^{+}_u|>1
$$
(see Appendix). The corresponding BVP is
\begin{equation}
\label{eqn:EF_original}
\left\{
      \begin{array}{rcl}
        \dot{w} + f_u^{\rm T}(x^{+},\alpha) w & = & 0 \;,\\
        w(T^{+}) - \mu w(0) & = & 0\;,\\
        \langle w(0),w(0)  \rangle - 1& = &0 \;,
      \end{array}
    \right.
\end{equation}
where $x^{+}$ is the solution of (\ref{eqn:L}).
In our implementation the above BVP is replaced by 
an equivalent BVP
\begin{equation}
\label{eqn:EF}
\left\{
      \begin{array}{rcl}
        \dot{w} +  f_u^{\rm T}(x^{+},\alpha) w + \lambda w& = & 0 \;,\\
        w(T^{+}) - s w(0) & = & 0\;,\\
        \langle w(0), w(0)  \rangle - 1& = & 0 \;,
      \end{array}
    \right.
\end{equation}
where $s={\rm sign\;}\mu=\pm 1$ and 
$$
\lambda=\ln|\mu|
$$ (see Appendix). In (\ref{eqn:EF}), the boundary conditions become
periodic or anti-periodic, depending on the sign of the multiplier
$\mu$, while the logarithm of its absolute value appears in the
variational equation. This ensures high numerical robustness.

Given $w$ satisfying (\ref{eqn:EF}), the projection boundary condition
(\ref{eqn:BVPp2cst}) becomes
\begin{equation}
\label{eqn:PROJCYCLE}
\langle w(0), u(\tau_{+}) - x^{+}(0) \rangle  = 0.
\end{equation}

\subsection{The connection}
Finally, we need a phase condition to select a unique periodic
solution among those which satisfy (\ref{eqn:L}), {\it i.e.}, to fix a
{\em base point} $x_0^+=x^+(0)$ on the cycle $O^{+}$ (see Figure
\ref{fig:conBVP}).  Usually, an integral condition is used to fix the
phase of the periodic solution.  For the point-to-cycle connection,
however, we need a new condition, since the end point near the cycle
should vary freely. To this end we require the end point of the
connection to belong to a plane orthogonal to the vector
$f_0^+=f(x^+(0),\alpha)$. This gives the following BVP for the
connecting solution:
\begin{equation}
\label{eqn:C}
    \left\{
      \begin{array}{rcl}
        \dot{u} - f(u,\alpha)  &=&  0 \;,\\
        \langle f(x^{+}(0),\alpha), u(\tau_{+}) - x^{+}(0) \rangle  &=&  0 \;.
      \end{array}
    \right.
\end{equation}

\subsection{The complete BVP}
\label{Sec:CompleteBVP}
The complete truncated BVP to be solved numerically consists of
(\ref{eqn:E1}), with 
\begin{equation}
\label{BVP-0}
u(0)=\xi + \varepsilon v,
\end{equation}
or (\ref{eqn:E2}), with 
\begin{subequations}
\label{BVP-1}
\begin{align}
\langle v, u(0)-\xi \rangle &= 0\;,\label{BVP-1-1}\\
\langle \eta, u(0)-\xi \rangle &= 0,\; \label{BVP-1-2}
\end{align}
\end{subequations}
as well as 
\begin{subequations}
\label{eqn:BVP}
\begin{align}
\dot{x}^+ - T^{+} f(x^+,\alpha) & = 0, \label{BVP-2}\\
x^+(0) - x^+(1) & = 0, \label{BVP-3}\\
\langle w(0), u(1) - x^{+}(0) \rangle  &= 0, \label{BVP-4}\\
\dot{w} + T^{+} f_u^{\rm T}(x^{+},\alpha) w + \lambda w & =  0, \label{BVP-5}\\
w(1) - s w(0) & =  0, \label{BVP-6}\\
\langle w(0), w(0)  \rangle - 1& =  0, \label{BVP-7}\\
\dot{u} - T f(u,\alpha) & = 0, \label{BVP-8}\\
\langle f(x^{+}(0),\alpha), u(1) - x^{+}(0) \rangle & =  0. \label{BVP-9}
\end{align}
\end{subequations}
Here the time variable is scaled to the unit interval $[0,1]$, so that
both the cycle period $T^{+}$ and the connecting time $T$ become
parameters.

If the connection time $T$ is fixed at a large value, this BVP allows
to continue simultaneously the equilibrium $\xi$, its eigenvalue
$\lambda_u$ or $\lambda_s$, the corresponding eigenvector $v$, the
periodic solution $x^+$ corresponding to the limit cycle $O^+$, its
period $T^{+}$, the logarithm of the absolute value of the unstable
multiplier of this cycle, the corresponding scaled eigenfunction $w$,
as well as (a truncation of) the connecting orbit $u$. These objects
become functions of one system parameter (when $\dim W^u_{-}=2$) or
two system parameters (when $\dim W^u_{-}=1)$.  These free system
parameters are denoted as $\alpha_i$.

If $\dim W^u_{-}=2$ then, generically, limit points (folds) are
encountered along the solution family. These can be detected, located
accurately, and subsequently continued in two system parameters, say,
$(\alpha_1,\alpha_2)$, using the standard fold-following facilities of
{\sc auto}.

\section{Starting strategies}\label{sec:homotopy}
The BVP's specified above can only be used if good starting data are
available. This can be problematic, since global objects -- a saddle
cycle and a connecting orbit -- are involved.  However, a series of
successive continuations in {\sc auto} can be used to generate all
necessary starting data, given little {\it a priori} knowledge about
the existence and location of a heteroclinic point-to-cycle
connection.

\subsection{The equilibrium and the cycle}
The equilibrium $\xi$, its unstable or stable eigenvalue, as well as
the corresponding eigenvector or adjoint eigenvector can be calculated
using {\sc maple} or {\sc matlab}. Alternatively, this saddle
equilibrium can often be obtained via continuation of a stable
equilibrium family through a limit point (fold) bifurcation.

To obtain the limit cycle $O^{+}$, one can continue numerically (with
\textsc{auto}\nocite{Doedeletal1997} or
\textsc{content}\nocite{KuznetsovLevitin1997}, for example) a limit
cycle born at a Hopf bifurcation to an appropriate value of $\alpha$,
from where we start the successive continuation.

\subsection{Eigenfunctions}\label{sec:eigfunc}
In the first of such continuations, the periodic solution
corresponding the limit cycle at the particular parameter values is
used to get an eigen\-function. To explain the idea, let us begin with
the original adjoint eigen\-function $w$.  Consider the periodic BVP
(\ref{BVP-2})--(\ref{BVP-3}) for the cycle, to which the standard
integral phase condition is added,
\begin{equation}
\label{eqn:CYCLEPHASE}
\int_0^1 \langle \dot{x}^{+}_{old}(\tau), x^{+}(\tau) \rangle = 0\;,
\end{equation}
as well as a BVP similar to (\ref{eqn:EF_original}), namely:
\begin{equation}
\label{eqn:EF_original_h}
\left\{
      \begin{array}{rcl}
        \dot{w} + T^{+} f_u^{\rm T}(x^{+},\alpha) w & = & 0 \;,\\
        w(1) - \mu w(0) & = & 0\;,\\
        \langle w(0),w(0)  \rangle - h& = &0 \;.
      \end{array}
    \right.
\end{equation}
In (\ref{eqn:CYCLEPHASE}), $x^{+}_{old}$ is a reference periodic
solution, typically the one in the preceding continuation step.  The
parameter $h$ in (\ref{eqn:EF_original_h}) is a {\em homotopy
parameter}, that is set to zero initially. Then
(\ref{eqn:EF_original_h}) has a trivial solution
$$
w(t) \equiv 0,\ \ h=0,
$$ for any real $\mu$. This family of trivial solutions parametrized
by $\mu$ can be continued in {\sc auto} using a BVP consisting of
(\ref{eqn:L}) (with scaled time variable $t$), (\ref{eqn:CYCLEPHASE}),
and (\ref{eqn:EF_original_h}) with free parameters $(\mu,h)$ and fixed
$\alpha$. A Floquet multiplier of the adjoint system then corresponds
to a branch point at $\mu_1$ along this trivial solution family (see
Appendix). {\sc auto} can accurately locate such a point and switch to
the nontrivial branch that emanates from it. Continuing this secondary
family in $(\mu,h)$ until, say, the value $h=1$ is reached, gives a
nontrivial eigenfunction $w$ corresponding to the multiplier $\mu_1$.
Note that in this continuation the value of $\mu$ remains constant,
$\mu \equiv \mu_1$, up to numerical accuracy.

The same method is applicable to obtain a nontrivial scaled adjoint
eigenfunction. For this, the BVP
\begin{equation}
\label{eqn:EF_h}
\left\{
      \begin{array}{rcl}
        \dot{w} + T^{+} f_u^{\rm T}(x^{+},\alpha) w + \lambda w & = & 0 \;,\\
        w(1) - s w(0) & = & 0\;,\\
        \langle w(0),w(0)  \rangle - h& = &0 \;,
      \end{array}
    \right.
\end{equation}
where $s={\rm sign}(\mu)$, replaces (\ref{eqn:EF_original_h}). A
branch point at $\lambda_1$ then corresponds to the adjoint multiplier
$s{\rm e}^{\lambda_1}$.  Branch switching then gives the desired
eigendata.

\subsection{The connection}\label{sec:conn}
Sometimes, an approximation of the connecting orbit can be obtained by
time-integration of (\ref{eqn:ODE}) with a starting point satisfying
(\ref{eqn:PROJEXP1}) or (\ref{eqn:PROJEXP2}) and (\ref{eqn:PROJEXP3}).
These data (the periodic solution corresponding to the limit cycle,
its nontrivial eigen\-function, and the integrated connecting orbit)
must then be merged, using the same scaled time variable and mesh
points. This only works for non-stiff systems provided that the
connecting orbit and its corresponding parameter values are known {\it
a priori} with high accuracy, which is not the case for most models.

A practical remedy in most cases is to apply the method of {\em
successive continuation} first introduced by Doe\-del, Fried\-man and
Mon\-tei\-ro (1993) for point-to-point problems. This method does not
guarantee that a connection will be found but works well if we start
sufficiently close to a connection {\em in the parameter space}.  Here
we generalize this method to point-to-cycle connections.

We first consider the case $\dim W^u_{-}=1$. To start, we
introduce a BVP composed of (\ref{eqn:E1}), (\ref{BVP-0}), and a
modified version of (\ref{eqn:BVP}), namely:
\begin{subequations}
\label{eqn:BVP-homotopy}
\begin{align}
\dot{x}^+ - T^{+} f(x^+,\alpha) & = 0, \label{BVP-2-h}\\
x^+(0) - x^+(1) & = 0, \label{BVP-3-h}\\
\Psi[x^{+}] &= 0, \label{BVP-4-h}\\
\dot{w} + T^{+} f_u^{\rm T}(x^{+},\alpha) w + \lambda w & = 0,
\label{BVP-5-h}\\
w(1) - s w(0) & =  0, \label{BVP-6-h}\\
\langle w(0), w(0)  \rangle - 1& =  0, \label{BVP-7-h}\\
\dot{u} - T f(u,\alpha) & = 0, \label{BVP-8-h}\\
\langle f(x^{+}(0),\alpha), u(1) - x^{+}(0) \rangle - h_1 & = 0,
\label{BVP-9-h}
\end{align}
\end{subequations}
where $\Psi$ in (\ref{BVP-4-h}) defines any phase condition fixing the
base point $x^{+}(0)$ on the cycle $O^+$; for example
$$
\Psi[x^{+}]=x^+_j(0)-a_j,
$$ where $a_j$ is the $j$th-coordinate of the base point at some given
parameter values, and $h_1$ is a {\em homotopy parameter}.

Take an initial solution to this BVP that collects the previously
found equilibrium-related da\-ta, the cycle-related da\-ta
$(x^{+},T^{+})$ including $x^{+}(0)$, the eigen\-function-related
da\-ta $(w,\lambda)$, as well as the value of $h_{1}$ computed for the
initial ``connection"
\begin{equation}
\label{eqn:INIT1}
u(\tau)=\xi + \varepsilon v {\rm e}^{\lambda_u T\tau},~~ \tau \in [0,1],
\end{equation}
which is a solution of the scaled linear approximation of
(\ref{eqn:ODE}) in the tangent line to the unstable manifold $W^u_{-}$
of $\xi$.  By continuation in $(T,h_1)$ for a fixed value of $\alpha$,
we try to make $h_1=0$, while $u(1)$ is near the cycle $O^+$, so that
$T$ becomes sufficiently large.

After this is accomplished, we introduce another BVP
composed of (\ref{eqn:E1}), (\ref{BVP-0}), and 
\begin{subequations}
\label{eqn:BVP-homotopy1}
\begin{align}
\dot{x}^+ - T^{+} f(x^+,\alpha) & = 0, \label{BVP-2-h1}\\
x^+(0) - x^+(1) & = 0, \label{BVP-3-h1}\\
\langle w(0), u(1) - x^{+}(0) \rangle  - h_2 &= 0, \label{BVP-4-h1}\\
\dot{w} + T^{+} f_u^{\rm T}(x^{+},\alpha) w + \lambda w & = 0,
\label{BVP-5-h1}\\
w(1) - s w(0) & =  0, \label{BVP-6-h1}\\
\langle w(0), w(0)  \rangle - 1& =  0, \label{BVP-7-h1}\\
\dot{u} - T f(u,\alpha) & = 0, \label{BVP-8-h1}\\
\langle f(x^{+}(0),\alpha), u(1) - x^{+}(0) \rangle & =  0, \label{BVP-9-h1}
\end{align}
\end{subequations}
where $h_2$ is another homotopy parameter.

Using the solution obtained in the previous step, we can activate one
of the system parameters, say $\alpha_{1}$, and aim to find a solution
with $h_2=0$ by continuation in $(\alpha_{1},h_2)$ for fixed $T$. Then
we can improve the connection by continuation in $(\alpha_{1},T)$,
restarting from this latest solution, in the direction of increasing
$T$.  Eventually, we fix a sufficiently large value of $T$ and
continue the (approximate) connecting orbit in two systems parameters,
say $(\alpha_{1},\alpha_{2})$, using the original BVP without any
homotopy parameter as described in Section \ref{Sec:CompleteBVP}. All
these steps are illustrated for the Lorenz example in Section
\ref{sec:Lorenz}. In practice, intermediate continuations in
$\varepsilon$ or other system parameters may be necessary to obtain a
good approximation to the connecting orbit.

When $\dim W^u_{-}=2$, a minor modification of the above homotopy
method is required. In this case, we replace (\ref{BVP-1}) by the
explicit boundary conditions
\begin{subequations}
\label{eqn:LEFT-h}
\begin{align}
u(0) - \xi - \varepsilon (c_1v^{(1)}+c_2v^{(2)}) &= 0, \label{LEFT-h-1}\\
c_1^2 +c_2^2 &= 1, \label{LEFT-h-2}
\end{align}
\end{subequations}
where $\varepsilon$ is a small parameter specifying the distance
between $u(0)$ and $\xi$, $v^{(j)}$ are two linear-independent vectors
tangent to $W^u_{-}$ of the saddle $\xi$, and $c_{1,2}$ are two new
scalar homotopy parameters. Note that if $v=(v_1,v_2,v_3)^{\rm T}$ is
a solution to (\ref{eqn:E2}) with $v_2 \neq 0$, then one can use the
normalized vectors
$$
v^{(1)}=\left(\begin{array}{r}
v_2 \\ -v_1 \\ 0~\end{array}\right), ~ 
v^{(2)}=\left(\begin{array}{r}
0~ \\ v_3 \\ -v_2\end{array}\right) .
$$ Now consider a BVP composed of (\ref{eqn:E2}), (\ref{eqn:LEFT-h}),
and (\ref{eqn:BVP-homotopy}). The initial data for this BVP are the
same as in the case $\dim W^u_{-}=1$, except for
$$
c_1=1,~~c_2=0.
$$
The initial ``connection" in this case is
\begin{equation}
\label{eqn:INIT2}
u(\tau)=\xi + \varepsilon {\rm e}^{\tau TA}v^{(1)} ,~~ \tau \in [0,1],
\end{equation}
where $A=f_u(\xi,\alpha)$, to be used to compute the initial value of
$h_{1}$ in (\ref{BVP-9-h}).

By continuation in $(T,h_1)$ (and, eventually, in $(c_1,c_2,h_1)$) for
fixed values of all other parameters, we aim to locate a solution with
$h_1=0$, with $u(1)$ near the base point of the cycle $O^+$, so that
$T$ becomes sufficiently large.  We then switch to the BVP composed of
(\ref{eqn:E1}), (\ref{eqn:LEFT-h}), and (\ref{eqn:BVP-homotopy1}), and
we aim to locate a solution with $h_2=0$, by continuation in
$(c_{1},c_{2},h_2)$ for fixed $T$. When this is achieved, we have a
solution to the original BVP (\ref{eqn:E2}), (\ref{BVP-1}), and
(\ref{eqn:BVP}) introduced in Section \ref{Sec:CompleteBVP} and
containing no homotopy parameters. Using this BVP, we can continue the
approximate connecting orbit in one system parameter, say
$\alpha_{1}$, with $T$ fixed.

Examples of such successive continuations will be given in Section
\ref{sec:3dfoodchain}, where we consider the standard model of a
3-level food chain. In that section also an alternative BVP
formulation for (\ref{eqn:LEFT-h}) is given. When one system parameter
is varied, limit points (folds) can be found and then continued in two
system parameters.

\section{Implementation in AUTO}\label{sec:implementation}

Our algorithms have been implemented in \textsc{auto}, which solves
the boundary value problems using superconvergent {\em orthogonal
collocation} with adaptive meshes.  {\sc auto} can compute paths of
solutions to boundary value problems with integral constraints and
non-separated boundary conditions:
\begin{subequations}
\label{eqn:AUTO}
\begin{align}
\dot{U}(\tau) - F(U(\tau),\beta) &=0\; , \; \; \tau \in [0,1], \label{5.1} \\
b(U(0),U(1),\beta) &= 0\;,  \label{5.2} \\
\int^1_0 q(U(\tau),\beta) d\tau &= 0\;, \label{5.3} 
\end{align}
\end{subequations}
where
$$
U(\cdot),F(\cdot,\cdot) \in  {\mathbb R}^{n_d}, ~
b(\cdot,\cdot) \in {\mathbb R}^{n_{bc}}, ~
q(\cdot,\cdot) \in {\mathbb R}^{n_{ic}}, 
$$
and
$$
\beta \in {\mathbb R}^{n_{fp}}. 
$$ Here $\beta$ represents the $n_{fp}$ {\it free parameters} that are
allowed to vary, where
\begin{equation}
n_{fp} = n_{bc} + n_{ic}  - n_d + 1.
\label{5.4}
\end{equation}
The function $q$ can also depend on $\dot{U}$ and on the derivative of
$U$ with respect to pseudo-arclength, as well as on $\hat{U}$, the
value of $U$ at the previously computed point on the solution family.

For our primary BVP problem (\ref{eqn:E1}) or (\ref{eqn:E2}) with
(\ref{BVP-0}) or (\ref{BVP-1}), respectively, and (\ref{eqn:BVP}), we
have
$$
n_d = 9, ~~ n_{ic}=0,
$$ and $n_{bc}=19$ or $18$, respectively, since (\ref{eqn:E1}) and
(\ref{eqn:E2}) are treated as boundary conditions.

\section{Examples}\label{sec:results}

In this section we illustrate the performance of our algorithm by
applying it to three model systems, namely, the Lorenz equations, an
electronic circuit model, and a biologically relevant system.

\subsection{The Lorenz system}\label{sec:Lorenz}
One of the best-known dynamical systems that has a heteroclinic
point-to-cycle connection is the three-dimensional Lorenz system,
given by
\begin{equation}\label{eqn:Lorenz}
\left\{
  \begin{array}{rcl}
    \dot{x}_1 & = &\sigma (x_2 - x_1) ,\\
    \dot{x}_2 & = & r x_1 - x_2 - x_1 x_3 ,\\
    \dot{x}_3 & = & x_1 x_2 - b x_3 ,
  \end{array}
\right.
\end{equation}
with standard parameter values $\sigma=10$, $b=8/3$, and where $r$ is
the usual bifurcation parameter. With these parameter values, a
supercritical pitchfork bifurcation from the trivial equilibrium
occurs at $r = 1$, giving rise to two symmetric nontrivial equilibria.
At $r \approx 13.962$ there are two symmetry-related orbits of
infinite period that are homoclinic to the origin, and from which two
families of saddle cycles arise (together with a nontrivial hyperbolic
invariant set).  A subcritical Hopf bifurcation of nontrivial
equilibria takes place at $r_H \approx 24.7368$, where these two
cycles disappear.
\begin{figure*}[ht]
\begin{center}
\includegraphics[width=17.0cm]{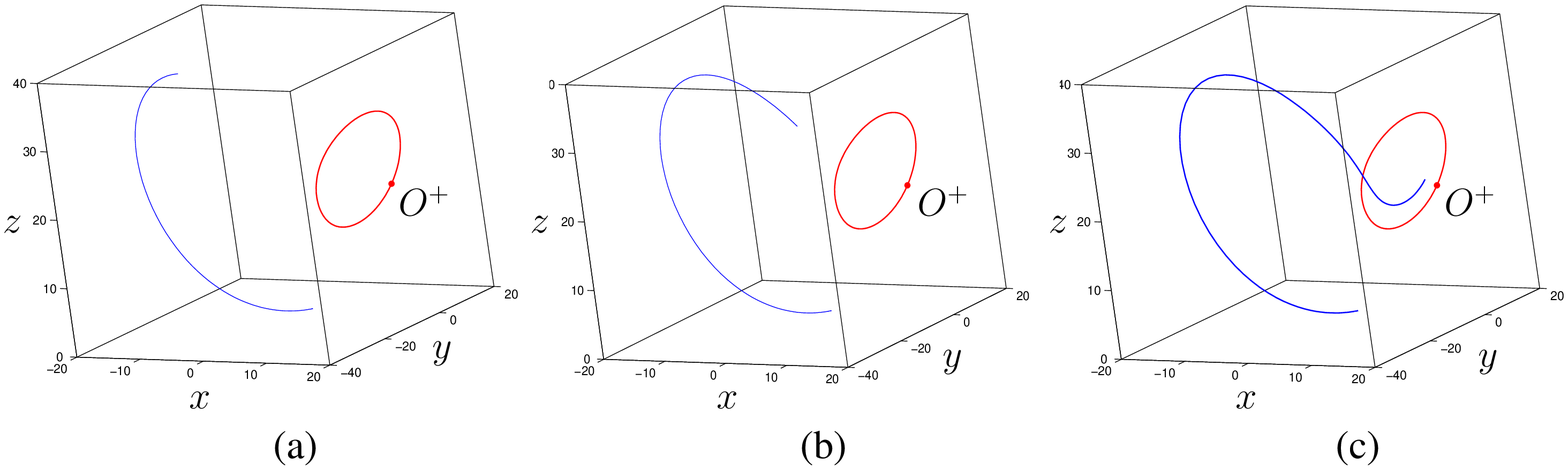}
\end{center}
\caption{Continuation in $T$: (a) $T=1.43924$; (b) $T=1.54543$; (c)
$T=2.00352$.}
\label{fig:LorHomT} 
\end{figure*}

At a critical value $r_{het}$ there is a heteroclinic point-to-cycle
connection, that generates a chaotic attractor, see
Afraimovich et al. (1977\nocite{AfrBykShil1977}). Its domain of
attraction is bounded by the stable invariant manifolds of the saddle
cycles.  Beyn (1990\nocite{Beyn1990}) found $r_{het} \approx 24.05$,
and later Die\-ci and Rebaza (2004\nocite{DieciRebaza2004}) calculated
$$ 
r_{het} = 24.057900322267 \ldots 
$$

The heteroclinic connection can be continued in two parameters, for
example $r$ and $\sigma$ with $b$ fixed.  The resulting curve in the
$r,\sigma$-plane was first shown in Appendix II, written by
L.P. Shil'\-ni\-kov, to the Russian translation of the book by
Mars\-den and Mc\-Crac\-ken (see Pampel (2001\nocite{Pampel2001}),
Die\-ci and Rebaza (2004), for more recent related results). 
As shown by Bykov and Shilnikov (1992\nocite{BykovShilnikov1992}),
the canonical Lorenz attractor appears by crossing only a part of the 
heteroclinic connection curve.

We begin at $r=21.0$ and consider a saddle limit cycle $O^{+}$ of
(\ref{eqn:Lorenz}) with the base point
$$
x^{+}(0)=(9.265335, 13.196014, 15.997250)
$$
and period $T^{+}=0.816222$. This cycle can be obtained easily by
continuation in {\sc auto} and has two nontrivial multipliers:
$$
\mu^{+}_s = 0.0000113431,~~\mu^{+}_u = 1.26094.
$$ 
\begin{figure}[htbp]
\includegraphics[width=9.5cm]{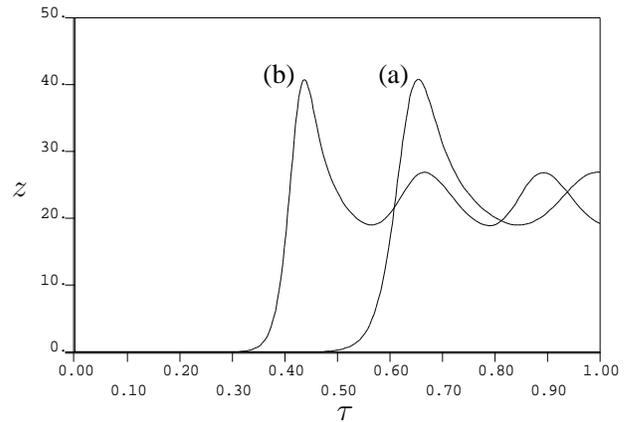}
\caption{Two profiles of the truncated connecting orbit in the Lorenz
  system scaled to the unit time interval: (a) $T=2.00352$; (b)
  $T=3.0$.}
\label{fig:connectlor} 
\end{figure}
\begin{figure*}[htbp]
\begin{center}
\includegraphics[width=12.5cm]{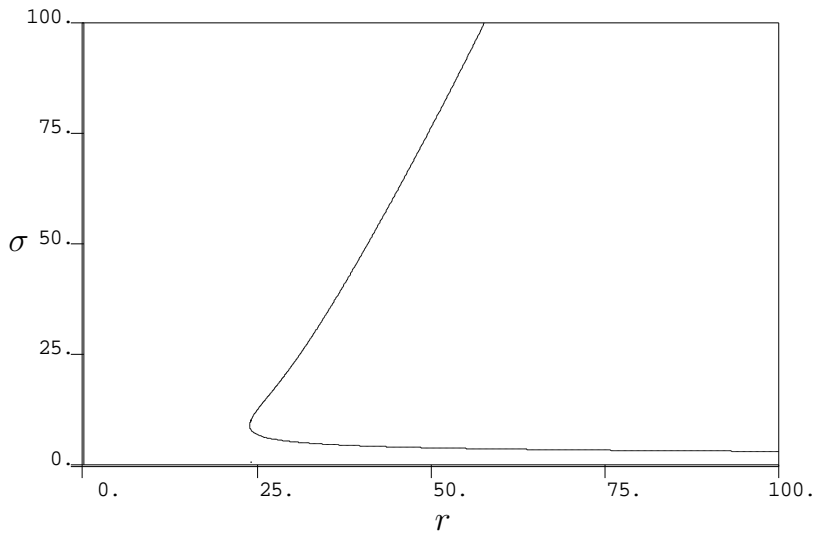}
\end{center}
\vspace{-1.0cm}
\caption{The bifurcation curve of the Lorenz system corresponding to
the point-to-cycle connection.}
\label{fig:parlor} 
\end{figure*}

To compute the eigenfunction $w$, we first continue the trivial
solution of the BVP (\ref{BVP-2}), (\ref{BVP-3}),
(\ref{eqn:CYCLEPHASE}), and (\ref{eqn:EF_h}), to detect a branch point
at
$$
\lambda=\ln(\mu^{+}_u)=0.231854,
$$ from which a nontrivial branch is followed until the value $h=1$ is
reached. This gives a nontrivial eigenfunction $w(t)$, with
$\|w(0)\|=1$, namely,
$$
w(0)=(0.168148, 0.877764, -0.448616)^{\rm T}.
$$ In these continuations all problem parameters, that is $r,\sigma,$
and $b$, are fixed.

The next step is to find an approximation to the connecting orbit. For
this, we consider the BVP (\ref{eqn:E1}), (\ref{BVP-0}), and
(\ref{eqn:BVP-homotopy}) with
$$
\Psi[x^{+}]=x_1^{+}(0) - 9.265335
$$ and continue its solution at fixed system parameters with respect
to $(T,h_1)$.  Figure \ref{fig:LorHomT} shows three consecutive
solutions with $h_1=0$. The end point of the last solution (with
$T=2.00352$) is located near the base point $x^{+}(0)$ of the cycle
$O^+$. Using this solution as the initial data for the BVP
(\ref{eqn:E1}), (\ref{BVP-0}), and (\ref{eqn:BVP-homotopy1}), we do a
continuation in $(r,h_2)$ with $T$ fixed until $h_2=0$ is detected.
This occurs at $r=24.0720$, and ensures that the end point of the
connection is in a plane orthogonal to $w(0)$, {\it i.e.}, in the
tangent plane to $W^s_{+}$ at $x^{+}(0)$.

The primary BVP consisting of (\ref{eqn:E1}), (\ref{BVP-0}), and
(\ref{eqn:BVP}) is used for further continuation runs.  First, the
length of the connecting orbit is increased by continuation in $(r,T)$
until $T=3.0$.  The corresponding parameter value $r=24.0579$ gives a
good approximation for $r_{het}$, since the `tail' of the connecting
orbit follows the cycle $O^+$ several times; (see Figure
\ref{fig:connectlor}).

Finally, continuation in the two system parameters $(r,\sigma)$ with
$T$ fixed, gives the bifurcation curve corresponding to the
point-to-cycle connection in (\ref{eqn:Lorenz}), see Figure
\ref{fig:parlor}.

\subsection{A circuit model}\label{sec:circuitmodel}

The next example is one from the \textsc{Homcont} demos of Champneys
{\it et al.} (1999\nocite{Champneysetal1996}), namely, the electronic
circuit model of Freire et al.  (1993\nocite{Freireetal1993}; see also
the \textsc{auto} demos \textit{tor} and \textit{cir}). The equations
are
\begin{equation}\label{eqn:Freireetal}
\left\{
  \begin{array}{rcl}
    r\dot{x}_1 & \!=\! &-(\beta\!+\!\nu) x_1 \!+\! \beta x_2
    \!-\! a_3 x_1^3 \!+\! b_3 (x_2\!-\!x_1)^3,\\
    \dot{x}_2 & \!=\! &\beta x_1 - (\beta \!+ \!\gamma)x_2
    - x_3 - b_3(x_2 \!-\! x_1)^3,\\
    \dot{x}_3 & \!=\! &x_2,
  \end{array}
\right.
\end{equation}
where $\gamma$ = 0, $r$ = 0.6, $a_3$ = 0.328578, $b_3$ = 0.933578, and
$\nu$ and $\beta$ are bifurcation parameters.  With \textsc{Homcont}
it was shown previously that a homoclinic connection to the origin
occurs for
$$
\nu_{init} = -0.721309 ~,~ \beta_{init} = 0.6
$$ with truncated time interval $T$ = 200. Continuation in
two-parameter dimension then leads to a Shil'nikov-Hopf bifurcation at
$$
\nu = -1.026445 ~,~ \beta = -2.330391 \cdot 10^{-5},
$$ where a limit cycle bifurcates from the equilibrium, effectively
turning the homoclinic connection into a heteroclinic one (see
\textsc{auto} demo \textit{cir}). We can now compare the results from
the continuation in \textsc{Homcont} with the results from the
application of our BVP system.

The equilibrium in this system is a saddle-focus, and we therefore
have $n_s^- = 2$ and $n_u^- = 1$. To generate appropriate starting
data we locate a Hopf bifurcation, with $\beta$ as free parameter,
from where a cycle is continued up to a selected value of $\beta$,
say, $\beta = -0.32$.  The saddle limit cycle $O^+$ has the base point
$$
x^+(0) = (0.03448278,0.46460323,0.4737975)
$$
and period $T^+ = 6.3646138$. The nontrivial multipliers are
$$
\mu^+_s = 3.986051\cdot10^{-6},\; \mu^+_u = 18.85438
$$

The eigenfunction of this cycle is computed as described in Section
\ref{sec:eigfunc}, which yields
$$
w(0) = (0.99950, -0.019205, 0.024767)^{\rm T}
$$
and the log multiplier
$$
\lambda = -13.579343187.
$$

An approximation of the connecting orbit is then obtained using BVP
(\ref{eqn:E1}), (\ref{BVP-0}), and (\ref{eqn:BVP-homotopy}), with
$$
\Psi[x^{+}]=x_2^{+}(0) - 0.46460323.
$$ The software \textsc{content} is used to get a good approximation
of the connection period $T$, after which shooting in \textsc{matlab}
is used to obtain the orbit itself for the given period.

Continuation of this approximate orbit with respect to $(T,h_1)$
yields several orbits with $h_1 = 0$. For $T = 11.59816$ the orbit is
close enough to the $x_2$ base coordinate to use the data for the BVP
(\ref{eqn:E1}), (\ref{BVP-0}), and (\ref{eqn:BVP-homotopy1}).
Continuation in $(\nu,h_2)$ is done until a zero of $h_2$ is reached.

The primary BVP (\ref{eqn:E1}), (\ref{BVP-0}), and (\ref{eqn:BVP}) is
used in the subsequent computations. Continuation in $(\nu,T)$ gives
orbits of any desired period $T$; we used $T = 20$ with
$$
\nu = -1.500498.
$$
At this point continuation can be done in $(\nu,\beta)$.
\begin{figure}
  \scalebox{0.5}{
    \includegraphics{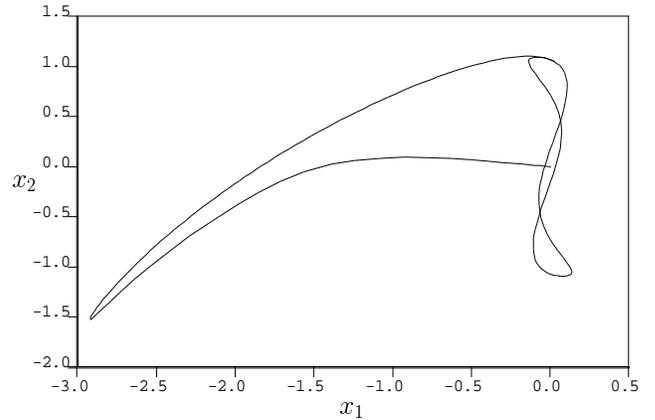}}
  \vspace{-1.0cm}
  \caption{A point-to-cycle connection of the electronic circuit
            model, projected onto the $x_1,x_2$-plane.}
  \label{fig:cirxy}
\end{figure}
\begin{figure}
  \scalebox{0.5}{
    \includegraphics{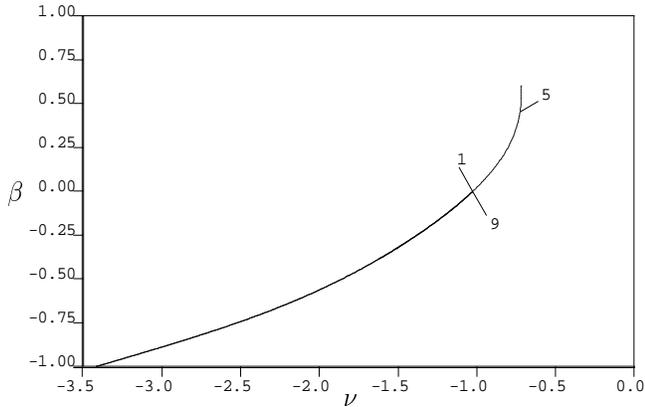}}
  \vspace{-1.0cm}
  \caption{Continuation in ($\nu$, $\beta$) of the point-to-cycle
           connection, as explained in detail in the text.}
  \label{fig:cirbif}
\end{figure}

In Figure \ref{fig:cirxy} we see a point-to-cycle connection in a
$x_1$,$x_2$-plot at some selected parameter values. It is apparent
that the homotopy method has resulted in a good approximation of the
connecting orbit.  Figure \ref{fig:cirbif} shows the composite results
of the two-parameter continuation of the homoclinic connection in
\textsc{Homcont} and our continuation of the heteroclinic connection.
Label 5 is the starting point of the continuation of the homoclinic
connection that terminates at the solution labelled~1. Beyond this
solution \textsc{Homcont} gives spurious results. Note that label 1
coincides with label 9, where the curve of the heteroclinic connection
turns back onto itself, {\it i.e.}, the continuation reverses
direction approximately at the point where the Shil'nikov-Hopf
bifurcation occurs. Plots in \textsc{auto} of the limit cycle data
(not shown) reveal that indeed the cycle shrinks practically to a
point, before the continuation reverses direction.

\subsection{A food chain model}\label{sec:3dfoodchain}

The following three-level food chain model from theoretical biology
is based on the Rosen\-zweig-MacArthur (1963\nocite{RosenzweigMacArthur1963})
prey-predator model. The equations are given by
\begin{equation}\label{eqn:RM3D}
  \left\{
    \begin{array}{rcl}
      \dot{x}_1 & = & x_1(1-x_1) - f_1(x_1,x_2),\\
      \dot{x}_2 & = & f_1(x_1,x_2) - f_2(x_2,x_3) - d_1 x_2,\\
      \dot{x}_3 & = & f_2(x_2,x_3) - d_2 x_3,
    \end{array}
  \right.
\end{equation}
with Holling Type-II functional responses
$$
f_i(u,v) = \frac{a_i u v}{1 + b_i u}\;,~~i=1,2.
$$

The death rates $d_1$ and $d_2$ are used as bifurcation parameters,
with the other parameters set at $a_1 = 5$, $a_2 = 0.1$, $b_1 = 3$,
and $b_2 = 2$.

It is well known that this model displays chao\-tic behaviour in a given
pa\-ra\-me\-ter range, see Ho\-ge\-weg and Hes\-per
(1978\nocite{HogewegHesper1978}), Kle\-ba\-noff and Has\-tings
(1994\nocite{KlebanoffHastings1994}), McCann and Yod\-zis
(1995\nocite{McCannYodzis1995}), Kuz\-net\-sov and Ri\-nal\-di
(1996\nocite{KuznetsovRinaldi1996}), and Kuz\-net\-sov et
al. (2001\nocite{Kuznetsovetal2001}).

Previous work by Boer et al.  (1999\nocite{Boeretal1999},
2001\nocite{Boeretal2001}) has also shown that the regions of chaos
are intersected by homoclinic and heteroclinic global connections.  
In particular, a heteroclinic point-to-cycle orbit connecting
a saddle with a two-dimensional unstable manifold to a saddle cycle
with a two-dimensional stable manifold can exist. It was shown that the stable manifold of
this limit cycle forms the basin boundary of the interior attractor
and that the boundary has a complicated structure, especially near the
equilibrium, when the heteroclinic orbit is present.  
These and other results were obtained
numerically using multiple shooting.  In this section we reproduce
these results for the heteroclinic point-to-cycle connection. Using
our homotopy method we obtain an accurate approximation of the
heteroclinic orbit.  A one-parameter bifurcation diagram then shows
limit points, which correspond to tangencies of the above-mentioned
two-dimensional manifolds. We then continue the limit points in two
parameters.

A starting point can be found, for example, at $d_1 \approx
0.2080452$, $d_2 = 0.0125$, where there is a fold bifurcation in which
two limit cycles appear.  This also corresponds to the birth of the
heteroclinic point-to-cycle connection.

Before using the homotopy method to obtain an approximation of the
point-to-cycle connection, we locate a Hopf bifurcation, for instance
at $d_1 \approx 0.51227$, $d_2 = 0.0125$.  The limit cycle born at
this Hopf bifurcation is continued up to a selected value of $d_1$,
say, $d_1 = 0.25$.

We now have an equilibrium 
$$
\xi =(0.74158162,   0.16666666,   11.997732 )
$$
and a saddle limit cycle with the base point
$$
x^+(0) = (0.839705,0.125349,10.55289)
$$
and period $T^+ = 24.282248$. Its nontrivial multipliers are
$$
\mu^+_s = 0.6440615, \mu^+_u = 6.107464\cdot10^{2}.
$$

The eigenfunction $w$ is obtained as described in the previous
sections. Continuation of the trivial solution of the BVP
(\ref{BVP-2}), (\ref{BVP-3}), (\ref{eqn:CYCLEPHASE}), and
(\ref{eqn:EF_h}) and the subsequent continuation of the bifurcating
family until $h = 1$, yields the multiplier
$$
\lambda = \ln(\mu^+_s) = -0.439961.
$$ Note that we use the stable multiplier, because of the projection
boundary conditions. The associated nontrivial eigenfunction $w(t)$
with $\|w(0)\| = 1$ has
$$
w(0) = (0.09306,-0.87791,-4.69689)^{\rm T}.
$$

We now consider a BVP composed of (\ref{eqn:E2}), (\ref{eqn:LEFT-h}a),
and (\ref{eqn:BVP-homotopy}).  Using \textsc{content} and
\textsc{matlab} we obtain an approximation of the connection with the
boundary condition
$$
\Psi[x^+] = x_2^+(0) - 0.125349
$$ and period $T = 155.905$. The starting point is calculated by
splitting the normalized adjoint stable vector (evaluated at $d_1 =
0.25,d_2 = 0.0125$)
$$
v = (0.098440,0.168771,0.0049532)^{\rm T}
$$ into $v^{(1)}$ and $v^{(2)}$, as described in Section
\ref{sec:conn}, and multiplying it by a small $\varepsilon$, say
$\varepsilon=0.001$. In our case the starting point was
$$
u(0) = (0.742445,0.166163,11.997732).
$$

The first homotopy step involves continuation in $(h_1,T)$. However,
this does not lead to zeroes of $h_1$. To obtain $h_1 = 0$ we expand
the previous set of BVPs with (\ref{eqn:LEFT-h}b). Subsequent
continuation in $(c_1,c_2,h_1)$ gives a solution with $h_1 = 0$ that
indeed ends near the base point $x^+(0)$ of the limit cycle.

For continuation in the second homotopy step, a switch is made to a
BVP composed of (\ref{eqn:E2}), (\ref{eqn:BVP-homotopy1}) and
(\ref{eqn:LEFT-h}).  Continuation in $(c_1,c_2,h_2)$ leads to some
solutions with $h_2 = 0$.

\begin{figure}
  \scalebox{0.5}{
    \includegraphics{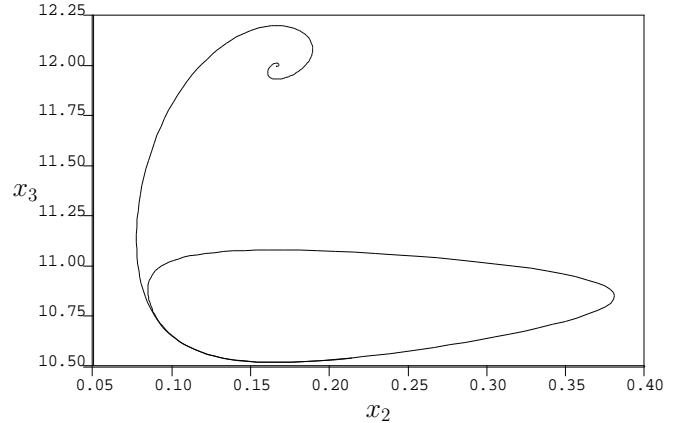}}
  \vspace{-1.0cm}
  \caption{An approximation to the point-to-cycle connection projected
    onto the $(x_2,x_3)$-plane for the food chain model
    with $a_1 = 5$, $a_2 = 0.1$, $b_1 = 3$, $b_2 = 2$, $d_1=0.25$, and
    $d_2=0.0125$.}
  \label{fig:connectrm200}
\end{figure}

\begin{figure}
  \scalebox{0.5}{
    \includegraphics{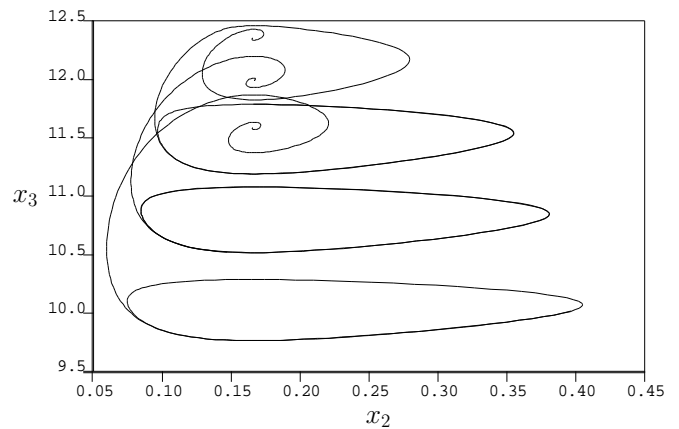}}
  \vspace{-1.0cm}
  \caption{Several point-to-cycle connections in the
  food chain model with different values of $d_1$.}
  \label{fig:connectrm201}
\end{figure}

The obtained approximate connecting point-to-cycle connection now
suffices for continuation in system parameters.  Before doing a
continuation in a system parameter the connection is improved by
increasing the connection period. A user-defined point of $T = 300$
suffices. Next, the parameter $\varepsilon$ is decreased up to a
user-defined point of $\varepsilon = -1 \cdot 10^{-5}$, so that the
starting point $u(0)$ is slightly away from the equilibrium $\xi$.
Figure~\ref{fig:connectrm200} displays a projection of the
point-to-cycle connection onto the ($x_2,x_3$)-plane.

Now the connecting orbit can be continued up to a limit point in one
system parameter. Figure~\ref{fig:connectrm201} displays three
connecting orbits obtained after continuation with respect to
$\alpha_1 = d_1$.  Continuations in $d_1$ result in the detection of
the points
$$
d_1 = 0.280913\quad \textrm{and} \quad d_1 = 0.208045
$$ where the first one is a limit point and the second one a
termination point. This point coincides with a tangent bifurcation for
the limit cycle to which the point-to-cycle orbit connects.
Continuations in $d_2$ result in the detection of the points
$$
d_2 = 0.0130272\quad \textrm{and} \quad d_2 = 9.51660 \cdot 10^{-3}
$$ which are both limit points.
\begin{figure*}
\begin{center}
\includegraphics[width=12.5cm]{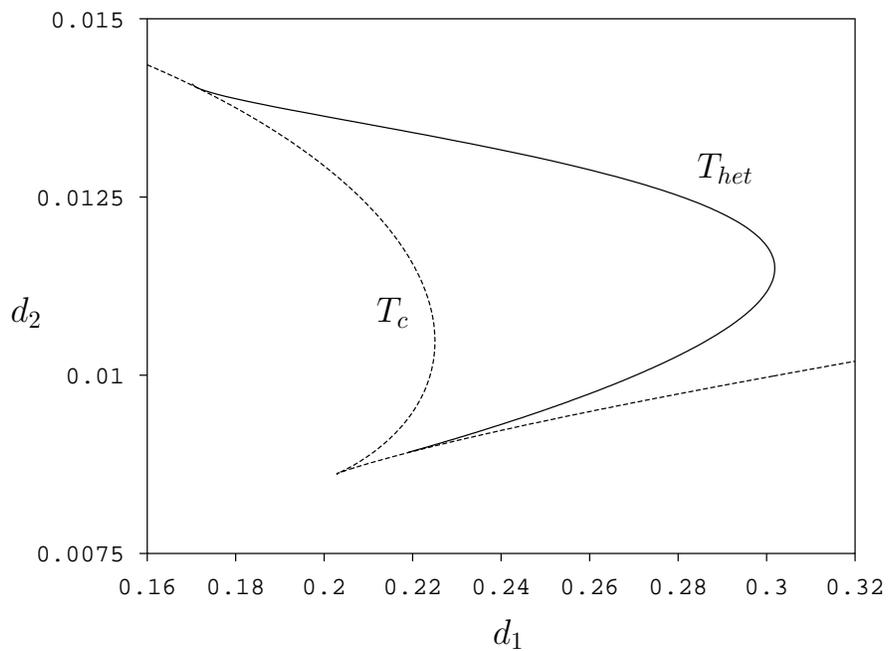}
\end{center}
  \vspace{-1.0cm}
  \caption{A two-parameter bifurcation diagram of the food chain model
    that shows the region where there exist point-to-cycle
    connections. The region is bounded on one side by the cycle fold,
    ($T_c$), and on the other side by the curve $T_{het}$, the locus
    of limit points of the heteroclinic connections.}
  \label{fig:hetcurve2d}
\end{figure*}
Any of the detected limit points can now be used as a starting point
for a two-parameter continuation in $\alpha = d_1$ or $d_2$. In practice,
the connection period may have to be increased or decreased to obtain
the full two-parameter continuation curve. In the demo, the last limit
point ($d_2 = 9.51660 \cdot 10^{-3}$) is the one selected for the food
chain model. The two-parame\-ter continuation curve terminates at both
ends in co\-dim 2 points lying on the above-mentioned tangent
bifurcation for the limit cycle. These points coincide with the log
multiplier $\lambda = 0$. Observe that this corresponds to the point
$d_1 = 0.208045$, detected in the one-parameter continuation, where
also $\lambda = 0$.

For the continuation in two system parameters, the BC
(\ref{eqn:LEFT-h}) proves ineffective, since it leads to the detection
of several spurious limit points. This is, because the orbit
spiralling out from the equilibrium has an elliptical shape. The
circle of a small radius, centered at the equilibrium, intersects the
spiral at several points, one of which is the starting point of the
connecting orbit. During continuation, with a changing problem
parameter, the spiral will change size and the starting point on the
circle may collide with another such point where the circle and the
spiral intersect. This intersection would correspond to a fold with
respect to the problem parameter. As a result, to obtain a full
continuation curve of the connecting orbit in two system parameters,
some restarts are required.

In order to avoid these spurious folds, we returned to the original 
BC (\ref{BVP-1}) with (\ref{BVP-1-2}) in the form
\begin{align}\label{eqn:bcdoed}
  u_j(0) - \xi_j  & = 0\;,
\end{align}
where $j$ is either 1,2 or 3. By setting $j = 2$ we are in line with the
work by Boer et al.\ (2001\nocite{Boeretal2001}), who used a
Poincar\'{e} plane through the equilibrium $\xi$ where $x_2 = \xi_2$.

Figure~\ref{fig:hetcurve2d} shows the curve of limit points $T_{het}$
that is computed with the method described above, using the standard
switching and fold-following facilities of \textsc{AUTO}. This curve
can be obtained in one run, given the connection period is chosen
conveniently. It agrees with the results previously obtained by Boer
et al.\ (1999\nocite{Boeretal1999}) by labourious multiple shooting.

\section{Discussion}

Our continuation method for point-to-cycle connections, using
homotopies in a boundary value setting is both robust and
time-efficient. Detailed {\sc auto} demos that carry out the
computations described in Section 6 are freely downloadable from \par
\noindent
{\tt www.bio.vu.nl/thb/research/project/globif}.  \par Although the
method was presented for 3D-systems, it can be extended directly to
point-to-cycle connections in $n$-dimensional systems, when the
unstable invariant manifold of the equilibrium $\xi$ is either
one-dimensional or has codimension one, while the stable invariant
manifold of the cycle $O^+$ has codimension one.

In the forthcoming Part II of this paper, we will extend our method to
include detection and continuation of cycle-to-cycle connections.

%% file: appendix.tex
\appendix

\section{Monodromy matrices}\label{app:monmatrix}

In order to approximate the invariant manifolds of a limit cycle we use eigenvalues
and eigenfunctions of appropriate variational problems. These eigenvalues in turn are
the eigenvalues of the so-called \textit{monodromy matrix}.

To define an eigenfunction of the periodic solution $x(t+T)$ = $x(t)$,
where $T$ is the period of the cycle, of an autonomous system of smooth ODE's
\begin{equation}
\label{eqn:ODE_0}
\dot{u}=f(u),\ \ \ f:\mathbb{R}^n \to \mathbb{R}^n,
\end{equation}
write a solution of this system near the cycle in the form
$$
u(t)  = x(t) + \xi(t)\;,
$$
where $\xi(t)$ is a small deviation from the periodic solution.
After substitution and truncation of the $O(\|\xi\|^2)$-terms, we obtain the following
{\em variational system}:
\begin{equation}
\label{eqn:varprob}
  \dot{\xi} = A(t) \xi\;\;,\quad \xi \in \mathbb{R}^n \;,
\end{equation}
where $A(t)=f_u(x(t))$ is the Jacobian matrix evaluated along the periodic solution;
$A(t+T)$ = $A(t)$.

Now, consider the matrix initial-value problem
\begin{equation}
\label{eqn:matfund}
  \dot{Y} = A(t) Y\;,\ \ \ \ Y(0)=I_n,
\end{equation}
where $I_n$ is the unit $n$ $\times$ $n$ matrix. Its solution $Y(t)$ at $t=T$ 
is the {\em monodromy matrix} of the cycle: 
$$
M=Y(T).
$$
The monodromy matrix is nonsingular.
Any solution $\xi(t)$ to (\ref{eqn:varprob}) satisfies
\begin{equation}\label{eqn:monodromy}
  \xi(T) = M \xi(0)\;.
\end{equation}
The eigenvalues of the
monodromy matrix $M$ are called the {\em Floquet multipliers} of the cycle.
There is always a multiplier $+1$. Moreover, the product of all multipliers
is positive:
$$
\mu_1 \mu_2 \cdots \mu_n=\exp\left(\int_0^{T} {\rm div\ }f(x(t))\ dt\right).
$$

Together with (\ref{eqn:varprob}), consider the {\em adjoint variational system}
\begin{equation}
\label{eqn:adjvarprob}
  \dot{\zeta} = -A^{\rm T}(t) \zeta\;\;,\quad \zeta \in \mathbb{R}^n \;
\end{equation}
and the corresponding matrix initial-value problem
\begin{equation}
\label{eqn:matfundadj}
  \dot{Z} = -A^{\rm T}(t) Z\;, \ \ \ \ Z(0)=I_n,
\end{equation}
which is the adjoint system to (\ref{eqn:matfund}). Note, that the
multipliers of the adjoint monodromy matrix 
$$
N=Z(T)
$$ 
are the {\em inverse} multipliers of the monodromy matrix $M=Y(T)$.
The proof of this well-known fact goes as follows. Compute
\begin{eqnarray*}
    \frac{d}{dt} (Z^{\rm T} Y)& = &\frac{d Z^{\rm ^T}}{dt} Y
    + Z^{\rm T}\frac{d Y}{dt}\\
    & = &(-A^{\rm T}Z)^{\rm T} Y + Z^{\rm T} A Y\\
    & = & Z^{\rm T} (-A) Y + Z^{\rm T} A Y = 0\;.
\end{eqnarray*}
Since  $Z(0) = Y(0) = I_n$, we get $Z^{\rm T}(T) Y(T)  = I_n$,
which implies
$$
N=[M^{-1}]^{\rm T}.
$$

Due to (\ref{eqn:monodromy}), a multiplier $\mu$ satisfies $v(T)=\mu v(0)$ with
$v(0)\neq 0$ or, equivalently, it is a solution component of
the following BVP on the unit interval $[0,1]$:
\begin{equation}
\label{eqn:EF0}
\left\{
      \begin{array}{rcl}
        \dot{v} + T A(t) v & = & 0 \;,\\
        v(1) - \mu v(0) & = & 0\;,\\
        \langle v(0),v(0)  \rangle - 1& = &0 \;.
      \end{array}
    \right.
\end{equation}
First assume that $\mu>0$ and write 
$$
\mu={\rm e}^{\lambda},\ \ v(t)={\rm e}^{\lambda t}w(t).
$$
Then $w$ satisfies a periodic BVP, namely:
\begin{equation}
\label{eqn:EF1}
\left\{
      \begin{array}{rcl}
        \dot{w} + T A(t) w + \lambda w & = & 0 \;,\\
        w(1) - w(0) & = & 0\;,\\
        \langle w(0),w(0)  \rangle - 1& = &0 \;.
      \end{array}
    \right.
\end{equation}
Similarly, when $\mu<0$, we can introduce
$$
\mu=-{\rm e}^{\lambda},\ \ v(t)={\rm e}^{\lambda t}w(t)
$$
and obtain an anti-periodic BVP
\begin{equation}
\label{eqn:EF2}
\left\{
      \begin{array}{rcl}
        \dot{w} + T A(t) w + \lambda w & = & 0 \;,\\
        w(1) + w(0) & = & 0\;,\\
        \langle w(0),w(0)  \rangle - 1& = &0 \;.
      \end{array}
    \right.
\end{equation}
This technique can easily be adapted to the multipliers of
the adjoint variational problem (\ref{eqn:adjvarprob}).

Finally, we note that the eigenvalue problem for a Floquet multiplier
$$
M v - \mu v  = 0
$$
can be considered as a {\em continuation problem} with $n+1$ variables $(v,\mu) \in \mathbb{R}^n
\times \mathbb{R}$ defined by $n$ equations. This continuation problem has a trivial
solution family $(v,\mu)=(0,\mu)$. An eigenvalue $\mu_1$ corresponds to a {\em branch point},
from which a {\em secondary solution family} $(v,\mu_1)$ with $v \neq 0$ emanates. This nontrivial
family can be continued using an extended continuation problem
$$
\left\{\begin{array}{rcl}
M v - \mu v  &=& 0\;,\\
\langle v,v \rangle - h &=&0,
\end{array}
\right.
$$
which consists of $n+1$ equation with $n+2$ variables $(v,\mu,h)$. If $h=1$ is reached,
we get a normalized eigenvector $v$ corresponding to the eigenvalue $\mu_1$, since
along this branch $\mu \equiv \mu_1$. Generalization of this procedure to the BVP
(\ref{eqn:EF0}) (as well as to (\ref{eqn:EF1}), (\ref{eqn:EF2}), and their adjoint versions)
is straightforward.

%% file: HET.bbl
\begin{thebibliography}{}

\bibitem[\protect\citename{Afraimovich \bgroup \em et al.\egroup ,
  }1977]{AfrBykShil1977}
V.~S. Afraimovich, V.~V. Bykov, and L.~P. Shilnikov, [1977], ``The origin and
  structure of the {L}orenz attractor,'' {\em Dokl. Akad. Nauk SSSR}, {\bf
  234}, 336--339.

\bibitem[\protect\citename{Beyn, }1990]{Beyn1990}
W.-J. Beyn, [1990], ``The numerical computation of connecting orbits in
  dynamical systems,'' {\em IMA J. Numer. Anal.}, {\bf 10}, 379--405.

\bibitem[\protect\citename{Beyn, }1994]{Beyn1994}
W.-J. Beyn, [1994], ``On well-posed problems for connecting orbits in dynamical
  systems.'', In {\em Chaotic {N}umerics (Geelong, 1993)}, volume 172 of {\em
  Contemp. Math.},  131--168. Amer. Math. Soc., Providence, RI.

\bibitem[\protect\citename{Boer \bgroup \em et al.\egroup ,
  }1999]{Boeretal1999}
M.~P. Boer, B.~W. Kooi, and S.~A. L.~M. Kooijman, [1999], ``Homoclinic and
  heteroclinic orbits to a cycle in a tri-trophic food chain,'' {\em J. Math.
  Biol.}, {\bf 39}, 19--38.

\bibitem[\protect\citename{Boer \bgroup \em et al.\egroup ,
  }2001]{Boeretal2001}
M.~P. Boer, B.~W. Kooi, and S.~A. L.~M. Kooijman, [2001], ``Multiple attractors
  and boundary crises in a tri-trophic food chain,'' {\em Math. Biosci.}, {\bf
  169}, 109--128.

\bibitem[\protect\citename{Bykov and Shilnikov, }1992]{BykovShilnikov1992}
V.~V. Bykov and A.~L. Shilnikov, [1992], ``On the boundaries of the domain of
  existence of the {L}orenz attractor,'' {\em Selecta Mathematica Sovietica},
  {\bf 11}, 375--382.

\bibitem[\protect\citename{Champneys and Kuznetsov,
  }1994]{ChampneysKuznetsov1994}
A.~R. Champneys and {\mbox{Yu. A}}.~Kuznetsov, [1994], ``Numerical detection
  and continuation of codimension-two homoclinic bifurcations,'' {\em Internat.
  J. Bifur. Chaos Appl. Sci. Engrg.}, {\bf 4}, 785--822.

\bibitem[\protect\citename{Champneys \bgroup \em et al.\egroup ,
  }1996]{Champneysetal1996}
A.~R. Champneys, {\mbox{Yu. A}}.~Kuznetsov, and B.~Sandstede, [1996], ``A
  numerical toolbox for homoclinic bifurcation analysis,'' {\em Internat. J.
  Bifur. Chaos Appl. Sci. Engrg.}, {\bf 6}, 867--887.

\bibitem[\protect\citename{Dieci and Rebaza, }2004a]{DieciRebaza2004a}
L.~Dieci and J.~Rebaza, [2004], ``Erratum: ``{P}oint-to-periodic and
  periodic-to-periodic connections'','' {\em BIT Numerical Mathematics}, {\bf
  44}, 617--618.

\bibitem[\protect\citename{Dieci and Rebaza, }2004b]{DieciRebaza2004}
L.~Dieci and J.~Rebaza, [2004], ``Point-to-periodic and periodic-to-periodic
  connections,'' {\em BIT Numerical Mathematics}, {\bf 44}, 41--62.

\bibitem[\protect\citename{Doedel and Friedman, }1989]{DoedelFriedman1989}
E.~J. Doedel and M.~J. Friedman, [1989], ``Numerical computation of
  heteroclinic orbits,'' {\em J. Comput. Appl. Math.}, {\bf 26}, 155--170.

\bibitem[\protect\citename{Doedel \bgroup \em et al.\egroup ,
  }1994]{DoFrMo1994}
E.~J. Doedel, M.~J. Friedman, and A.~C. Monteiro, [1994], ``On locating
  connecting orbits,'' {\em Appl. Math. Comput.}, {\bf 65}, 231--239.

\bibitem[\protect\citename{Doedel \bgroup \em et al.\egroup ,
  }1997]{Doedeletal1997}
E.~J. Doedel, A.~R. Champneys, T.~F. Fairgrieve, {\mbox{Yu. A}}.~Kuznetsov,
  B.~Sandstede, and X.~Wang, [1997], ``{\sc auto97}: {C}ontinuation and
  bifurcation software for ordinary differential equations.'', Technical
  report, Concordia University, Montreal, Quebec, Canada.

\bibitem[\protect\citename{Freire \bgroup \em et al.\egroup ,
  }1993]{Freireetal1993}
E.~Freire, A.~J. Rodr{\'{\i}}guez-Luis, E.~Gamero, and E.~Ponce, [1993], ``A
  case study for homoclinic chaos in an autonomous electronic circuit. {A} trip
  from {T}akens-{B}ogdanov to {H}opf-{\v{s}}il'nikov,'' {\em Phys. D}, {\bf
  62}, 230--253.

\bibitem[\protect\citename{Hogeweg and Hesper, }1978]{HogewegHesper1978}
P.~Hogeweg and B.~Hesper, [1978], ``Interactive instruction on population
  interactions.,'' {\em Computational Biology and Medicine}, {\bf 8}, 319--327.

\bibitem[\protect\citename{Klebanoff and Hastings,
  }1994]{KlebanoffHastings1994}
A.~Klebanoff and A.~Hastings, [1994], ``Chaos in three-species food chains,''
  {\em J. Math. Biol.}, {\bf 32}, 427--451.

\bibitem[\protect\citename{Kuznetsov and Levitin, }1997]{KuznetsovLevitin1997}
{\mbox{Yu. A}}.~Kuznetsov and V.~V. Levitin, [1997], ``{\sc CONTENT}:
  {I}ntegrated environment for the analysis of dynamical systems.''.
\newblock Centrum voor Wiskunde en Informatica (CWI), Kruislaan 413, 1098 SJ
  Amsterdam, The Netherlands.
\newblock Available for download at {\tt ftp://ftp.cwi.nl/pub/content}.

\bibitem[\protect\citename{Kuznetsov and Rinaldi, }1996]{KuznetsovRinaldi1996}
{\mbox{Yu. A}}.~Kuznetsov and S.~Rinaldi, [1996], ``Remarks on food chain
  dynamics,'' {\em Math. Biosci.}, {\bf 134}, 1--33.

\bibitem[\protect\citename{Kuznetsov \bgroup \em et al.\egroup ,
  }2001]{Kuznetsovetal2001}
{\mbox{Yu. A}}.~Kuznetsov, O.~De~Feo, and S.~Rinaldi, [2001], ``Belyakov
  homoclinic bifurcations in a tritrophic food chain model.,'' {\em SIAM J.
  Appl. Math.}, {\bf 62}, 462--487.

\bibitem[\protect\citename{McCann and Yodzis, }1995]{McCannYodzis1995}
K.~McCann and P.~Yodzis, [1995], ``Bifurcation structure of a three-species
  food chain model.,'' {\em Theor. Pop. Biol.}, {\bf 48}, 93--125.

\bibitem[\protect\citename{Pampel, }2001]{Pampel2001}
T.~Pampel, [2001], ``Numerical approximation of connecting orbits with
  asymptotic rate,'' {\em Numer. Math.}, {\bf 90}, 309--348.

\bibitem[\protect\citename{Rosenzweig and MacArthur,
  }1963]{RosenzweigMacArthur1963}
M.~L. Rosenzweig and R.~H. MacArthur, [1963], ``Graphical representation and
  stability conditions of predator-prey interactions.,'' {\em Am. Nat.}, {\bf
  97}, 209--223.

\end{thebibliography}
